\documentclass[12pt]{article} 
\usepackage{index}
\usepackage[reqno]{amsmath}
\usepackage{fixmath}
\usepackage{amssymb, amsthm, enumerate}
\usepackage{mathtools}

\usepackage{hyperref}

\DeclarePairedDelimiter{\abs}{\lvert}{\rvert}
\DeclarePairedDelimiter{\norm}{\lVert}{\rVert}

\newtheoremstyle{plainsl}%
	{\topsep}
	{\topsep}
	{\slshape} 
	{}
	{\normalfont\bfseries}
	{.}
	{ }
	{}

\swapnumbers

{\theoremstyle{plainsl}
\newtheorem{theorem}{Theorem}[section]
\newtheorem{lemma}[theorem]{Lemma}
\newtheorem{corollary}[theorem]{Corollary}}
{\theoremstyle{remark}
}

\newcommand\cref[1]{Corollary~\ref{cor:#1}}

\renewcommand\proof{\noindent\textsl{Proof. }}
\newcommand\sqr[2]{{\vbox{\hrule height.#2pt
    \hbox{\vrule width.#2pt height#1pt \kern#1pt
        \vrule width.#2pt}\hrule height.#2pt}}}
\renewcommand\qed{%
	\ifmmode\eqno\sqr53
	\else\nolinebreak\ \hfill\sqr53\medbreak\fi}

\numberwithin{equation}{section}

\newcommand\De{\Delta}

\newcommand\ga{\gamma}

\renewcommand\th{\theta} 

\newcommand\cx{{\mathbb C}}

\newcommand\rats{{\mathbb Q}}
\newcommand\comp[1]{{\mkern2mu\overline{\mkern-2mu#1}}}
\newcommand\diff{\mathbin{\mkern-1.5mu\setminus\mkern-1.5mu}}

\newcommand\seq[3]{#1_{#2},\ldots,#1_{#3}}

\newcommand\pmat[1]{\begin{pmatrix} #1 \end{pmatrix}}

\DeclareMathOperator{\tr}{tr}

\DeclareMathOperator{\sign}{sign}

\newcommand\one{{\bf1}}

\newcommand\hm{\widehat M}

\newcommand\psd{\succcurlyeq}

\newcommand\bA{\comp{A}}

\title{Sedentary Quantum Walks} 

\author{
	Chris Godsil\footnote{
		Research supported by Natural Sciences and Engineering Council of Canada, 
		Grant No. RGPIN-9439}\\
	Combinatorics \& Optimization\\
	University of Waterloo}

\begin{document}
\maketitle
	
\begin{abstract}
    Let $X$ be a graph with adjacency matrix $A$. The \textsl{continuous quantum walk}
	on $X$ is determined by the unitary matrices $U(t)=\exp(itA)$. If $X$ is the complete
	graph $K_n$ and $a\in V(X)$, then 
	\[
		1-|U(t)_{a,a}|\le2/n.
	\]
	In a sense, this means that a quantum walk on a complete graph stay home with 
	high probability. In this paper we consider quantum walks on cones over an $\ell$-regular
	graph on $n$ vertices. We prove that if $\ell^2/n\to\infty$ as $n$ increases, than a 
	quantum walk that starts on the apex of the cone will remain on it with probability 
	tending to $1$ as $n$ increases. On the other hand, if $\ell\le2$ we prove that
	there is a time $t$ such that local uniform mixing occurs, i.e., all vertices
	are equally likely. 
	
	We investigate when a quantum walk on strongly regular graph has a high probability
	of ``staying at home'', producing large families of examples with the stay-at-home 
	property where the valency is small compared to the number of vertices.
\end{abstract}

\section{Introduction}

Let $X$ be a graph with adjacency matrix $A$. A \textsl{continuous quantum walk} on $X$
is defined by the 1-parameter family of unitary matrices
\[
	U(t) = \exp(itA).
\]
We view these matrices as operating on the Hilbert space $\cx^{V(X)}$. If the state of the
corresponding system at time $0$ is given by a unit vector $z$, then the state at time $t$
is $U(t)z$. 

We use $M\circ N$ to denote the Schur product of matrices $M$ and $N$ (of the same order).
The questions concerning a quantum walk of physical interest can always be expressed as questions
concerning the \textsl{mixing matrices}
\[
	M(t) := U(t)\circ\comp{U(t)} = U(t)\circ U(-t)
\]
The entries of $M(t)$ are non-negative and each row and each column sums to 1. Thus
they specify probability densities on $V(X)$, and these describe the outcome of measurements
on a associated quantum system.

If $X$ is the complete graph on $n$  vertices, then continuous quantum walks on $X$
have a very surprising property: for large $n$, for any time $t$, the mixing matrix
$M(t)$ is close to $I$. Informally we say that $K_n$ has the ``stay-at-home-property''.

In this paper we do two things. A graph $X$ on $n+1$ vertices is 
a \textsl{cone} over the graph $Y$ on $n$ vertices if $Y$ is an induced subgraph of $X$ and
the vertex in $V(X)\diff V(Y)$ has valency $n$, that is, it is adjacent to each vertex in $Y$.
Thus $K_{1,n}$ is the cone over the empty graph on $n$ vertices and $K_n$ is the cone over $K_{n-1}$.
We investigate when a cone over a regular graph has the stay at home property.

The second question is which strongly regular graphs have the stay-at-home property.
A graph is strongly regular if it is regular (but neither complete nor empty) and
there are integers $a$ and $c$ such that each pair of adjacent vertices has exactly
$a$ common neighbours and each pair of distinct non-adjacent vertices has exactly
$c$ common neighbours. The graphs $mK_n$ formed by $m$ vertex disjoint copies
of $K_n$ is strongly regular whem $m,n>1$; more interesting examples are $C_5$
and the Petersen graph. Since strongly regular graph are the most regular graphs we have
(complete and empty graphs aside), considering them in this context is very natural.

\section{Complete Graphs}

We establish the stay-at-home property for the complete graphs. The result is not new,
but the machinery we develop to prove the result will be useful.

Recall that if $A$ is real symmetric matric with distinct eigenvalues 
\[
	\seq\th1m,
\]
then there are spectral projections $\seq E1m$ such that
\[
	A = \sum_r \th_r E_r.
\]
Hece $E_r$ represents orthogonal projection onto the $\th_r$-eigenspace. Hence
\[
	E_r^T = E_r = E_r^2
\]
and $E_rE_s=0$ if $r\ne s$. Further $I=\sum_r E_r$ and, more generally, if $f$
is a function defined on the eigenvalues of $A$, then
\[
	f(A) = \sum_r f(\th_r) E_r.
\]

The eigenvalues of $K_n$ are $n-1$ (with multiplicity $1$) and $-1$ (with multiplicity $n-$);
the corresponding spectral idempotents are
\[
	E_0 = \frac1n J,\quad E_1 = I- \frac1n J.
\]
Therefore
\[
	U(t) = e^{(n-1)it}E_0 + e^{-it}E_1 = e^{-it}\Bigl(I +\frac{e^{nit}-1}{n}J\Bigr)
\]
and it follows immediately that the entries of $U(t)-e^{-it}I$ are bounded above in absolute
value by $\frac{2}n$, while the diagonal entries are bounded below in absolute value by $1-\frac{2}n$.

Note that if we can establish that there is a bound $c/n$
for all off-diagonal entries of $U(t)$ then, since $M(t)$ is doubly stochastic, we have
\[
	\abs{M(t)_{a,a}} \ge 1 -(n-1)\frac{c^2}{n^2} = 1 - \frac{c^2}{n} + \frac1{n^2},
\]
whence it follows that $X$ is stay-at-home. We will make free use of this observation.

\begin{lemma}
	\label{lem:comp-stay}
	If $X$ is a regular graph on $n$ vertices, then for any two vertices $a$ and $b$,
	\[
		\left\lvert\left(U_{\comp{X}}(t)-e^{-it}U_X(t))_{a,b}\right)\right\rvert \le \frac2n.
	\]
\end{lemma}

\proof
Assume $X$ has valency $k$. Since $\bA=J-I-A$, we see that
\[
	U_{\comp{X}}(t) = \exp(-itA)\exp(itI) \exp(itJ) = e^{-it}U_X(-t)\exp(itJ)
\]
Using the spectral decomposition of $J$ we have
\[
	\exp(itJ) = e^{nit} \frac1n J + I-\frac1n J = \frac{e^{nit}-1}n J + I.
\]
and since
\[
	U_X(-t) J = \sum_m \frac{(it)^m}{m!} (-A)^mJ = e^{-ikt}J,
\]
we conclude that
\[
	U_X(-t)\exp(itJ) = \frac{e^{nit}-1}n e^{-ikt} J + U_X(-t).\qed
\]

This result implies that if $X$ is regular, then it is stay-at-home if and only if $\comp{X}$ is.
In particular it yields another proof that $K_n$ has the stay-at-home property.

\section{Eigenvalues and Eigenvectors of Joins}
\label{sec:evevj}

Let $X$ be a $k$-regular graph on $m$ vertices and let $Y$ be an $\ell$-regular graph
on $n$ vertices. In this section we describe the spectral decomposition of their join.
Note that the join has an equitable partition $\pi$ with cells $(V(X),V(Y))$. 
Set $A=A(X)$ and $B=A(Y)$ and let $\hat{A}$ denote the adjacency matrix of $Z$. 
If $Z:=X+Y$ then the adjacency matrix for $Z$ is
\[
    \hat{A} = \pmat{A&J\\ J^T&B}
\]
and the adjacency matrix of the quotient $Z/\pi$ is
\[
	Q = \pmat{k&n\\ m&\ell}.
\]
Its eigenvalues are the zeros of the quadratic
\[
	t^2 - (k+\ell)t + (k\ell-mn),
\]
thus they are
\[
	\frac12(k+\ell \pm \sqrt{(k-\ell)^2+4mn}).
\]
We denote them by $\mu_1$ and $\mu_2$, with $\mu_1>\mu_2$.

Since
\[
	(Q-\mu_1I)(Q-\mu_2I) = 0 
\]
we see that the columns of $Q-\mu_2I$ are eigenvectors for $Q$ with eigenvalue $\mu_1$,
and the columns of $Q-\mu_1I$ are eigenvectors for $Q$ with eigenvalue $\mu_2$.
Hence the eigenvectors of $Z$ belonging to $\mu_1$ and $\mu_2$ respectively can be written
in partitioned form:
\[
	\pmat{(k-\mu_2)\one \\ m\one},\quad \pmat{(k-\mu_1)\one \\ m\one}.
\]
The remaining eigenvectors of $Z$ can be taken to be orthogonal to these two vectors,
and therefore such eigenvectors must sum to zero on $V(X)$ and $V(Y)$. If $x$ is an eigenvector
for $X$ orthogonal to $\one$ with eigenvalue $\lambda$, then
\[
	\pmat{x\\0}
\]
is an eigenvector for $Z$ with eigenvalue $\lambda$. Similarly if $y$ is an eigenvector
for $Y$ orthogonal to $\one$, then
\[
	\pmat{0\\y}
\]
is an eigenvector for $Z$ (with the same eigenvalue as $y$),

\section{Spectral Idempotents for Joins}

We are going to construct a refinement of the spectral decomposition of the join $Z$ of
$X$ and $Y$. (If $X$ and $Y$ are connected and have no eigenvalue in common, this will
be the actual spectral decomposition of $Z$.) This decomposition will, in large part,
be built from the spectral decompositions of $X$ and $Y$.

If $X$ is connected, we will use the spectral decomposition of $A$:
\[
	A = \sum_r \theta_r E_r
\]
where $\theta_r=k$ and $E_r=\frac1m J$. If $X$ is not connected, then $k$ has multiplicity
greater than one, and we may decompose the eigenspace belonging to $k$ as the sum of the span
of the constant vectors and the span of the vectors that are constant on components and sum 
to zero. The idempotent belonging to $k$ can then be written as the sum of $\frac1m J$
and a second idempotent. Now we have a refinement of the spectral decomposition of $X$,
with one extra term. We will still write this in the form above, with the understanding that 
$\theta_2=\theta_1$. A similar fuss can be made if $Y$ is not connected; we write its decomposition
as
\[
	B = \sum_s \nu_s F_s.
\]

If $Az=\th z$ and $\one^T z=0$, then
\[
    \hat{A}\pmat{z\\0} = \th\pmat{z\\0}.
\]
Similarly if $Bz=\th z$ and $\one^Tz=0$, then
\[
    \hat{A}\pmat{0\\z} = \th\pmat{0\\z}.
\]
We see that $n+m-2$ of the eigenvalues of $X+Y$ are eigenvalues of $X$ and eigenvalues
of $Y$.

Define
\[
    \hat{E}_r = \pmat{E_r&0\\0&0},\quad \hat{F}_s = \pmat{0&0\\0&F_s}
\]
and let $N_1$ and $N_2$ be the projections belonging to the eigenvalues $\mu_1$ and $\mu_2$
of $Z$. Then we have a decomposition for $Z$:
\begin{equation}
	\label{eq:hAmumu}
	\hat{A} =\mu_1 N_1 +\mu_2 N_2 +\sum_{r>1}\th_r\hat{E}_r +\sum_{s>1}\nu_s\hat{F}_s.
\end{equation}
We determined $\mu_1$ and $\mu_2$ in the previous section.
 
Since $\sum_r E_r=I$ and $\sum_s F_s=I$ we have
\[
	\sum_{r>1}\hat{E}_r = \pmat{I-\frac1m J&0\\0&0},\quad
	\sum_{s>1}\hat{F}_s = \pmat{0&0\\0&I-\frac1n J}
\]
and since the sum of the idempotents in \eqref{eq:hAmumu} is $I$, it follows that
\[
	N_1+N_2 = I - \pmat{I-\frac1m J&0\\ 0&I-\frac1n J}.
\]

The idempotent $N_1$ represents projection onto the span of the eigenvector
\[
	\pmat{(k-\mu_2)\one \\ m\one}
\]
and consequently
\[
	N_1  = c \pmat{ (k-\mu_2)^2 J_{m,m}& m(k-\mu_2)J_{m,n}\\ 
					 m(k-\mu_2)J_{n,m}& m^2 J_{n,n}}
\]
where $c$ is determined by the constraint $\tr(N_1)=1$. This means that
\[
	c^{-1} = m(k-\mu_2)^2 + m^2n = m((k-\mu_2)^2+mn).
\]
If we set 
\[
	\Delta = (k-\ell)^2+4mn
\]
then, after some calculation, we find that
\[
	(k-\mu_2)^2+mn = \sqrt{\Delta}(k-\mu_2).
\]
Hence $c^{-1}=m\sqrt{\Delta}(k-\mu_2)$. We can carry out similar calculations for $N_2$,
with the result that
\begin{align*}
	N_1  &= \frac1{m\sqrt{\Delta}(k-\mu_2)} \pmat{ (k-\mu_2)^2 J_{m,m}& m(k-\mu_2)J_{m,n}\\ 
					 m(k-\mu_2)J_{n,m}& m^2 J_{n,n}},\\[2pt]
	N_2  &= \frac1{m\sqrt{\Delta}(\mu_1-k)} \pmat{ (k-\mu_1)^2 J_{m,m}& m(k-\mu_1)J_{m,n}\\ 
					 m(k-\mu_1)J_{n,m}& m^2 J_{n,n}}.
\end{align*}

\section{The Transition Matrix of a Join}
\label{sec:tmatjoin}

Suppose $Z$ is the join $X+Y$ and $a$ and $b$ are two vertices in $X$. We want to determine
when we have perfect state transfer from $a$ to $b$ in $Z$. We note that if we do have
perfect state transfer from $a$ to $b$ at time $t$, then $U_Z(t)_{a,y}=0$ for all vertices
$u$ of $Y$ and $U_Z(t)_{a,a}=0$.

\begin{lemma}
	\label{lem:joinay}
	Assume $Z$ is the join of graphs $X$ and $Y$. If $a,b\in V(X)$ and $y\in V(Y)$, then
	\[
		U_{X+Y}(t)_{a,y} = \frac1{\sqrt{\Delta}}(\exp(i\mu_1 t) - \exp(i\mu_2 t))
	\]
	and
	\[
		U_{X+Y}(t)_{a,b} - U_X(t)_{a,b}
			= \frac1{m}\left( \frac{k-\mu_2}{\sqrt{\Delta}}\exp(i\mu_1 t) 
				- \frac{k-\mu_1}{\sqrt{\Delta}}\exp(i\mu_2 t) - \exp(ikt)\right).
	\]
\end{lemma}

\proof
Since $(\hat{E}_r)_{a,y}=(\hat{F_r})_{a,y}=0$ we have
\[
	U_{X+Y}(t)_{a,y} = \exp(i\mu_1 t)(N_1)_{a,y} + \exp(i\mu_2 t)(N_2)_{a,y}
		= \frac1{\sqrt{\Delta}}(\exp(i\mu_1 t) - \exp(i\mu_2 t)).
\]
From our spectral decomposition,
\[
	U_{X+Y}(t)_{a,b} = \frac{k-\mu_2}{m\sqrt{\Delta}}\exp(i\mu_1 t) - \frac{k-\mu_1}{m\sqrt{\Delta}}\exp(i\mu_2 t)
		+ \sum_{r>1}(E_r)_{a,b} \exp(i\theta_r t).
\]
Since $X$ is regular,
\[
	U_X(t) = \frac1m \exp(ikt)J_n + \sum_{r>1}\exp(i\theta_r t) E_r,
\]
from which our second expression follows.\qed

Therefore $U_{X+Y}(t)_{a,y}=0$ if and only if $\exp(it(\mu_1-\mu_2))=1$, that is,
if and only if for some integer $c$,
\[
	t = \frac{2c\pi}{\mu_1-\mu_2} = \frac{2c\pi}{\sqrt{\Delta}}.
\]
Since $(k-\mu_2)-(k-\mu_1)=\sqrt{\Delta}$, we find that for these values of $t$ we have
\[
	U_{X+Y}(t)_{a,b} - U_X(t)_{a,b} = \frac1m (\exp(i\mu_1t)-\exp(ikt)).
\]

\section{Irrational Periods and Phases}
\label{sec:irratpps}

We divert from our main theme to consider some consequences our the results in
the previous section. If $a$ and $b$ are distinct vertices in $X$, we have 
\textsl{perfect state transfer} from $a$ to $b$ if there is a time $t$ and a complex
scalar $\ga$ such that
\begin{equation}
	\label{eq:Utea}
	U(t)e_a = \ga e_b.
\end{equation}
The scalar $\ga$ is called a \textsl{phase factor}.
Since $U(t)$ is unitary, $\norm{\ga e_b}=1$ and therefore $\abs{\ga}=1$. If \eqref{eq:Utea}
holds, then
\[
	\ga^{-1}e_a = U(-t)e_b
\]
and, taking complex conjugates of both sides, we obtain
\[
	\ga e_a = U(t)e_b
\]
This show that if we have perfect state transfer from $a$ to $b$ at time $t$, we have
perfect state transfer from $b$ to $a$ at the same time, and with the same phase factor.
We also see that
\[
	U(2t)e_a = \ga^2 e_a.
\]
We say that $X$ is \textsl{periodic} at the vertex $a$ if there is a positive time $t$
such that $U(t)e_a$ is a scalar times $e_a$. Any vertex involved in perfect state transfer
is necessarily periodic.

In all known cases of perfect transfer, the phase factor $\ga$ is a root of unity.
We use the results of the Section~\ref{sec:tmatjoin} to construct examples of periodic vertices 
where the period is irrational and the associated phases are irrational.

A graph $X$ on $n$ vertices is a \textsl{cone} over a graph $Y$ if there is
a vertex $u$ of $X$ with degree $n-1$ such that $X\diff u\cong Y$. 
(Equivalently $X$ is isomorphic to $K_1+Y$.) We say
$u$ is the \textsl{apex} of the cone. 

\begin{lemma}
	Suppose $Y$ is an $\ell$-regular graph on $n$ vertices and let $Z$ be the cone over $Y$.
	The $Z$ is periodic at the apex with period $2\pi/\sqrt{\ell^2+4n}$.
\end{lemma}

\proof
We view $Z$ as the join of $X$ and $Y$, with $X=K_1$.
If $a$ is the apex vertex, then $\abs{U_Z(t)_{a,a}}=1$ if and only if
$U_Z(t)_{a,y}=0$ for all $y$ in $V(Y)$. This holds if and only if $t/2\pi$
is an integer multiple of $\Delta$, where $\Delta=\sqrt{\ell^2+4n}$.\qed

We have
\[
	U(t)_{a,a} = \sum_r e^{it\theta_r} (E_r)_{a,a}
\]
where $(E_r)_{a,a}\ge0$ and $\sum_r (E_r)_{a,a}=1$. Hence $|U(t)_{a,a}|=1$ if and only if
$e^{it\theta_r}=e^{it\theta_1}$ for all $r$, and therefore
\[
	U(t)_{a,a} = e^{it\theta_1}.
\]
So for a periodic cone the phase factor at the period is $e^{it\mu_1}$, where
\[
	\mu_1 = \frac12(\ell+\sqrt{\ell^2+4n}).
\]
Now it follows that for some integer $c$
\[
	t\mu_1 = \frac{2c\pi}{\sqrt{\ell^2+4n}} \frac12(\ell+\sqrt{\ell^2+4n})
		= \frac{c\pi\ell}{\sqrt{\ell^2+4n}} +c\pi.
\]
In all currently known cases where we have perfect state transfer, the phase factor is
a root of unity. The calculations we have just completed show that if $\ell^2+4n$
is not a perfect square, we have periodicity on the cone over $Y$ with phase factor not 
a root of unity.

\section{Walks on Cones}

We say a complex vector or matrix is \textsl{flat} if all its entries have the
same absolute value. We use $e_a$ for $a$ in $V(X)$ to denote the standard basis vectors
for $\cx^{V(X)}$. We have \textsl{uniform mixing relative to the vertex $a$} if there is a time $t$  
such that $U(t)e_a$ is flat. 
If we have uniform mixing relative to each vertex
at the same $t$, we say that $X$ admits \textsl{uniform mixing}. This holds if and
only if the matrix $U(t)$ is flat. Clearly uniform mixing is the antithesis of the stay-at-home
property.

There are many examples of graphs that admit uniform
mixing, including the hypercubes, but very few examples of graphs that admit uniform mixing
at a vertex at time $t$, but do not have uniform mixing at time $t$.
Carlson et al.~\cite{carlson-unimix} showed that there is uniform mixing on the star $K_{1,n}$, starting 
from the vertex of degree $n-1$.
We add to the list of examples. We say a graph $X$ on $n+1$ vertices is 
a \textsl{cone} over the graph $Y$ on $n$ vertices if $Y$ is an induced subgraph of $X$ and
the vertex in $V(X)\diff V(Y)$ has valency $n$, that is, it is adjacent to each vertex in $Y$.
Thus $K_{1,n}$ is the cone over the empty graph on $n$ vertices. 

\begin{lemma}
	If $Y$ is a regular graph with valency at most two and $Z$ is the cone over $Y$, then $Z$
	admits local uniform mixing starting from the apex.
\end{lemma}

\proof
Assume $n=|V(Y)|$ and that $Y$ is $\ell$-regular. Denote the cone over $Y$ by $Z$ and
let $a$ denote the apex. We set
\[
	\Delta = \sqrt{\ell^2+4n}
\]
and recall from Section~\ref{sec:evevj} that the eigenvalues in the eigenvalue support 
of $a$ are
\[
	\frac12(\ell \pm \Delta).
\]
We denote these by $\mu_1$ and $\mu_2$, assuming that $\mu_1>\mu_2$.

If $y\in V(Y)$ then Lemma~\ref{lem:joinay} yields that
\[
	U_Z(t)_{a,y} = \frac1\Delta \left(e^{it\mu_1}-e^{it\mu_2}\right)
		= \frac{e^{it\mu_2}}{\Delta} \left(e^{it\Delta}-1\right).
\]
We note that this is independent of the choice of $y$ in $Y$, and conclude that we have
uniform mixing from $a$ if and only if there is a time $t$ such that
\[
	\frac{1}{\Delta} \abs{e^{it\Delta}-1} = \frac1{\sqrt{n+1}};
\]
equivalently we need
\[
	\abs{e^{it\Delta}-1} = \frac{\sqrt{\ell^2+4n}}{\sqrt{n+1}}.
\]
As $\ell^2+4n=\ell^2-4 +4n+4$, the ratio on the right lies in the interval $[0,2]$
if and only if $\ell\le2$. Hence in these cases we can find a time $t$ for which satisfies
this equation, and then we have uniform mixing starting from $a$.\qed

By taking Cartesian powers, we get further examples of cones with uniform mixing starting 
from one vertex of the graph. 
It seems plausible that, in most cones, we do not get uniform mixing starting from a vertex in the 
base, but we do not have a proof in general. We leave it as an exercise to show that we do 
get uniform mixing in $K_{1,n}$ starting from a vertex of degree one if and only if $n=3$.
This was first observed by Hanmeng Zhan (private communication),
Cartesian powers of $K_{1,3}$ are our only known examples of graphs that admit uniform mixing
and are not regular.

\begin{corollary}
	Assume $Z$ is the cone over an $\ell$-regular graph $Y$ on $n$ vertices, with apex $a$.
	Then
	\[
		\abs{U(t)_{a,a}} \ge \frac{\ell^2}{\ell^2+4n}.
	\]
\end{corollary}

\proof
If $y\in V(Y)$ then from the proof of the theorem,
\[
	U_Z(t)_{a,y} = \frac{e^{it\mu_2}}{\Delta} \left(e^{it\Delta}-1\right).
\]
and consequently
\[
	\abs{U_Z(t)_{a,y}} \le \frac{2}{\sqrt{\ell^2+4n}}.
\]
Hence 
\[
	M(t)_{a,y} \le \frac{4}{\ell^2+4n}
\]
and (since $M(t)$ has row sum 1),
\[
	M(t)_{a,a} \ge 1 - \frac{4n}{\ell^2+4n} = \frac{\ell^2}{\ell^2+4n}.\qed
\]

Consider a sequence of $\ell$-regular graphs on $n$ vertices, where $\ell/n\to\infty$
as $n$ increases. Then, for large $n$, the corresponding cones have the stay at home
property at the apex.

We note also that if $\ell/\sqrt{n}\to c$ as $n\to\infty$, then $\abs{U(t)_{a,a}}\ge c/(c+4)$
while the off-diagonal entries $U(t)_{a,y}$ are bounded above by
\[
	\frac{2}{\sqrt{c+4}}n^{-1/2}.
\]
So a relaxed form of the stay-at-home property holds in this case.

\section{Mixing}

The \textsl{average mixing matrix} $\hm$ of a walk is defined by
\[	
	\hm = \lim_{T\to\infty}\frac1T \int_0^T M(t)\,dt.
\]
From \cite{cg-average} we have that
\[
	\hm = \sum_r E_r^{\circ2};
\]
it follows that $\hm$ is positive semidefinite and doubly stochastic.

For $K_n$ we have
\begin{align*}
	M(t) &= \Bigl(I +\frac{e^{nit}-1}{n}J\Bigr) \circ \Bigl(I +\frac{e^{-nit}-1}{n}J\Bigr) \\
		&= I + 2\frac{\cos(nt)-1}{n}I + \frac{2-2\cos(nt)}{n^2}J \\
		&= \Bigl(1-2\frac{1-\cos(nt)}{n}\Bigr)I +2\frac{1-\cos(nt)}{n^2}J
\end{align*}
while
\[
	\hm = \Bigl(1-\frac{2}{n}\Bigr)I + \frac{2}{n^2}J.
\]

\begin{lemma}
	\label{lem:im2mhi}
	We have $I\psd M(t) \psd 2\hm-I$;
\end{lemma}

\proof
First
\[
	M(t) = U(t)\circ U(-t) = \sum_{r,s} e^{it(\th_r-\th_s)} E_r\circ E_s
\]
and since $E_r\circ E_s=E_s\circ E_r$, it follows that
\begin{equation}
	\label{eq:Mtsum}
	M(t) = \sum_r E_r^{\circ2} + 2 \sum_{r<s}\cos((\th_r-\th_s)t).
\end{equation}
Now
\[
	I = \Bigl(\sum_r E_r\Bigr)^{\circ2} = \sum_{r,s} E_r\circ E_s
\]
and consequently
\begin{equation}
	\label{eq:imtsum}
	I-M(t) = \sum_{r<s}(2-2\cos((\th_r-\th_s)t))\, E_r\circ E_s
\end{equation}
The matrices $E_r\circ E_s$ are positive semidefinite and the above shows that $I-M$
ia a non-negative linear combination of positive semidefinite matrices, whence it
is positive semidefinite.

For the second inequality, we note first that
\[
	I = I\circ I = \sum_r E_r^2 + 2\sum_{r<s}E_r\circ E_s
\]
and, adding this to \eqref{eq:Mtsum} yields that
\begin{equation}
	\label{eq:mti2sum}
	M(t)+I = 2\sum_r E_r^2 + 2\sum_{r<s}(1+\cos((\th_r-\th_s)t)\,E_r\circ E_s.
\end{equation}
Appealing again to the fact that the Schur products $E_r\circ E_s$ are positive semidefinite,
we deduce that $M(t)-2\hm+I\psd0$.\qed

\begin{corollary}
	\label{cor:hmmt}
	If  $\hm_{a,a}=1-c$, then $M(t)_{a,a}\ge1-2c$ for all $t$.
\end{corollary}

\proof
If $M(t)-(2\hm-I)$ is positive semidefinite, then $M(t)_{a,a}\ge 2\hm_{a,a}-1$.\qed

If $M(t)>1-2c$ for all $t$, then clearly the diagonal entries of $\hm$ are at least $1-2c$.
So to prove that a graph is stay-at-home, it suffices to show that $\hm$ is close to $I$.

We consider what happens when one of the bounds in Lemma~\ref{lem:im2mhi} is tight.
A graph $X$ is periodic if $U_X(t)$ is a periodic function of $t$. Clearly $X$ is periodic
if its eigenvalues are integers; in particular $K_n$ is periodic.

\begin{lemma}
	If $X$ is connected and $I=M(t)$, then $X$ is periodic.
\end{lemma}

\proof
If $M(t)=I$, then $U(t)$ is diagonal. Since $U(t)$ commutes with $A$ and $X$ is
connected, $U(t)=\ga I$ for some $\ga$. Now $\det(U(t)=1$ because $\tr(A)=0$,
and therefore if $n=|V(X)|$, we have have $\ga^n=1$. Hence 
\[
	I = U(t)^n = U(nt)
\]
and therefore $X$ is periodic.\qed

\begin{lemma}
	If $X$ is connected and $M(t)=2\hm-I$, then $X$ is a complete graph.
\end{lemma}

\proof
Assuse $M(t)=2\hm-I$. By \eqref{eq:mti2sum} we have $\cos(t(\th_r-\th_s))=-1$, which
implies that $\cos(2t(\th_r-\th_s)t)=1$ and now \eqref{eq:imtsum} yields that $M(2t)=I$.

If $\cos(t(\th_r-\th_s))=-1$, then there is an odd integer $m_{r,s}$ such that
\[
	t(\th_r-\th_s) = m_{r,s}\pi
\]
and consequently
\[
	\frac{\th_r-\th_s}{\th_k-\th_\ell}
\]
is the ratio of two odd integers. Since 
\[
	\th_r-\th_s = (\th_r-\th_q) - (\th_q-\th_s),
\]
we see $X$ has at most two distinct eigenvalues; therefore it has diameter one and it
is a complete graph.\qed

From our expressions at the start of this session, if $X=K_n$ and $\cos(nt)=-1$ then
$M(t)=2\hm-I$.

\section{Absolute Value Bounds}

We derive an upper bound on the absolute values of the entries of $U(t)$.
If $a\in V(X)$, the \textsl{eigenvalue support} of $a$ is the set of eigenvalues
$\th_r$ such that $E_re_a\ne0$. Since the idempotents $E_r$ are positive semidefinite,
we have $E_re_a=0$ if and only if 
\[
	(E_r)_{a,a} = e_a^TE_re_a = 0.
\]
We say that a set $S$ of eigenvalues satisfies the \textsl{ratio condition} if
whenever 
\[
	\th_k,\ \th_\ell,\ \th_r,\ \th_s \in S
\]
and $\th_k\ne\th_\ell$, we have
\[
	\frac{\th_r-\th_s}{\th_k-\th_\ell} \in \rats.
\]

\begin{lemma}
	If $a,b\in V(X)$, then 
	\[
		\abs{U(t)_{a,b}} \le \sum_r (E_r)_{a,b}.
	\]
	If equality holds at time $t$, the intersection of eigenvalue supports of $a$ and $b$ 
	satisfies the ratio condition and $U(t)_{a,b}$ is a periodic function of $t$. 
\end{lemma}

\proof
Let $S$ denote the intersection of the eigenvalue supports of $a$ and $b$,

We have
\[
	U(t)_{a,b} = \sum_r e^{it\th_r} (E_r)_{a,b}
\]
and applying the triangle inequality, we find that
\[
	\abs{U(t)_{a,b}} \le \abs{(E_r)_{a,b}}.
\]
This is the stated bound. 

Equality holds in our bound if and only if the complex numbers
\[
	e^{it\th_r}\sign((E_r)_{a,b})
\]
are equal for all $r$ such that $(E_r)_{a,b}\ne0$. In this case,  $e^{2it(\th_r-\th_s)}=1$ 
for all eigenvalues $\th_r$ and $\th_s$ in $S$. Hence we have
\[
	U(2mt)_{a,b} = \sum_{r,s} (E_r)_{a,b} = 0
\]
for all $m$. We also see that there are integers $m_{r,s}$ such that
\[
	t(\th_r-\th_s) = m_{r,s}\pi
\]
and therefore, if $\th_k\ne\th_\ell$ and $\th_k$, $\th_\ell$, $\th_r$, $\th_s$
lie in $S$, then
\[
	\frac{\th_r-\th_s}{\th_k-\th_\ell} \in \rats.
\]

From the proof of\cite[Theorem~6.1]{cg-when} we have that either all elements of $S$ are 
integers, or they are all of the form $\frac12(a+b_i\sqrt\De$),
where $\De$ is a square-free integer and $a$ and $b_i$ are integers.\qed

\section{Strongly Regular Graphs: Spectral Decomposition}
\label{sec:srg-sd}

A graph $X$ is \textsl{strongly regular} if there are parameters $k$, $a$ and $c$ such that
\begin{enumerate}[(a)]
	\item 
	$X$ is $k$-regular.
	\item
	Each pair of adjacent vertices has exactly $a$ common neighbours.
	\item
	Each pair of distinct non-adjacent vertices have exactly $c$ common neighbours.
	\item
	$0 < k < n-1$.
\end{enumerate}
If $n=|V(X)|$ with $k$, $a$ and $c$ as given, we refer to the 4-tuple $(n,k;a,c)$
as the parameters of $X$. We use $\ell$ to denote $n-1-k$, the valency of the complement
of $X$.

In this section we review some of the properties of strongly regular graphs, including
the spectral decomposition of the adjacency matrix.

If $m,n>1$, then the disjoint union $mK_n$ of $m$ copies of $K_n$ is strongly regular.
The Petersen graph provides a more interestng example. The complement of a strongly
regular graph is strongly regular. A strongly regular graph is \textsl{primitive}
if it is not isomorphic to $mK_n$ or its complement. For more information of
strongly regular graphs see \cite{aeb-wh-spectra,cg-blue,cggfr-yellow}.

We denote the adjacency matrix of the complement of $X$ by $\bA$. A graph $X$ is
strongly regular with parameters $(n,k;a,c)$ if and only if
\[
	A^2 = kI +aA +c\bA
\]
or, equivalently
\[
	A^2- (a-c)A - (k-c)I = cJ,
\]
It follows that a primitive strongly regular graph has exactly three eigenvalue:
its valency $k$ and the two roots of the quadratic
\[
	t^2-(a-c)t-(k-c).
\]
These are
\[
	\frac12\left(a-c \pm\sqrt{(a-c)^2+4k-4c}\right)
\]
where the larger root is $\th$ and the smaller $\tau$. We denote the multiplicity of
$\th$ and $\tau$ respectively by $m_\th$ and $m_\tau$. These multiplcicities can
be expressed in terms of $n$, $k$, $\th$ and $\tau$:
\[
	m_\th = \frac{(n-1)(-\tau)-k}{\th-\tau}, \qquad m_\tau = \frac{(n-1)\th+k}{\th-\tau}.
\]

The adjacency matrix $A$ has spectral decomposition
\[
	A = kE_0 +\th E_1 +\tau E_2
\]
where $E_0 = \frac1nJ$,
\[
	E_1 = \frac{m_\th}n \left(I + \frac{\th}{k}A - \frac{\th+1}{\ell}\bA\right)
\]
and
\[
	E_2 = \frac{m_\tau}n \left(I + \frac{\tau}{k}A - \frac{\tau+1}{\ell}\bA\right).
\]

\section{Orthogonal Arrays, Steiner Designs}

With the exception of the so-called conference graphs, the eigenvalues
of strongly regular graphs are integers. Neumaier \cite{} proves a deep and remarkable 
result: for each positive integer $k$, there are only finitely many strongly regular graphs 
with least eigenvalue $-k$ that do not arise from either block graphs of Steiner $2$-designs 
or orthogonal arrays.

An \textsl{orthogonal array} $OA(k,n)$ is an $n^2\times k$ array with entries from $N=\{1,\ldots,n\}$
such that each ordered pair from $N\times N$ occurs exactly once in each pair of columns.
We define the graph of the array to be the graph with the rows of the array as its vertices,
with two rows adjacent if they agree on exactly one coordinate. The graph of an $OA(k,n)$
is strongly regular with parameters
\[
	n^2,\ k(n-1),\ n-2 +(k-1)(k-2),\ k(k-1).
\]
The eigenvalues are
\[
	k(n-1),\ n-k,\ -k
\]
with respective multiplicities
\[
	1,\ k(n-1),\ (n-1)(n+1-k).
\]
We note that $k\le n+1$ and if $k=n+1$, the graph is complete. We assume implicitly that $2\le k\le n$. 
For more information, see e.g \cite[Section~10.4]{cg-blue}.

\begin{lemma}
	\label{lem:oa-sed}
	If $X$ is the graph of an orthogonal array $OA(k,n)$ and $k=o(n)$, then $X$ has the 
	stay-at-home property.
\end{lemma}

\proof
The diagonal entries of $\hm$ are equal to
\[
	\frac{1}{n^2 } + \Bigl(\frac{k(n-1)}{n^2}\Bigr)^2 + \Bigl(\frac{(n-1)(n+1-k)}{n^2}\Bigr)^2.
\]
If $k=o(n)$, then $\hm\circ I$ converges to $I$ as $n$ increases, and Corollary~\ref{cor:hmmt}
now yields the conclusion.\qed

The parameters of the block graph of a Steiner $2$-$(v,k,1)$ design are
\[
	\frac{v(v-1)}{k(k-1)},\ k\frac{v-k}{k-1},\ \frac{v-1}{k-1}-2+(k-1)^2,\ k^2;
\]
its eigenvalues are
\[
	k\frac{v-k}{k-1},\ \frac{v-k}{k-1}-k,\ -k
\]
with respective multiplicities
\[
	1,\ v-1,\ \frac{v(v-1)}{k(k-1)} -v -1.
\]
Arguing as above, we conclude that $k=o(v)$, then the block graph is stay-at-home.

\section{Other Strongly Regular Graphs}

We show that, for a large strongly regular graph, the off-diagonal entries of $M(t)$
are small.

\begin{theorem}
	If $X$ is a strongly regular graph on $n$ vertices, the off-diagonal entries of $U(t)$
	are bounded in absolute value by $d/n^{1/4}$ for some $d$.
\end{theorem}

\proof
Assume $X$ is strongly regular. From Section~\ref{sec:srg-sd} we have
\[
	\frac{m_\th\th}{nk} = \frac{(n-1)(-\tau)\th-k\th)}{nk(\th-\tau)}
		= \frac{(n-1)(k-c)}{nk(\th-\tau)} - \frac{\th}{n(\th-\tau)}.
\]
Since $\tau<0$,
\[
	\frac{\th}{n(\th-\tau)} < \frac1n.
\]
while
\[
	\frac{(n-1)(k-c)}{nk(\th-\tau)} = \frac{n-1}{n} \frac{k-c}{k} \frac1{\th-\tau} <\frac1{\th-\tau}.
\]

From \cite[Section~6.1.8]{aeb-wh-spectra} or \cite[Lemma~10.3.1]{cggfr-yellow}
we have that
\[
	(\th-\tau)^2 = n\frac{k\ell}{m_\th m_\tau}.
\]
Since $X$ is $k$-regular with diameter two, $n\le k^2+1$. Therefore
\[
	(\th-\tau)^2 \ge n\frac{(n-\sqrt{n})\sqrt{n}}{n^2/4}
\]
and we conclude that, for some constant $d$,
\[
	\frac1{\th-\tau} < \frac{d}{n^{1/4}}.\qed
\]

We now consider the diagonal elements of $\hm$. From our expressions for the
spectral idempotents at the end of Section~\ref{sec:srg-sd}, it follows
that the diagonal entriers of $\hm$ are equal to
\[
	\frac1{n^2}(1 + m_\th^2 + m_\tau^2).
\]
Since $m_\th+m_\tau=n-1$, this is bounded below by
\[
	\frac1{n^2}(1 +\frac12(n-1)^2)
\]
which is less than $1/2$ (when $n\ge3$). In this case we get no udeful information
about the diagonalk entries of $M(t)$. We note that this lower bound is tioght
for conference graphs, where $m_\th=m_\tau=(n-1)/2$.

If $X$ is the graph of an $OA(k,n)$ and $k=\ga n$ with $0<\ga<1$, then from the expression
at the start of the proof of Lemma~\ref{lem:oa-sed}, we see that the diagonal entries
of $\hm$ tend to
\[
	\ga^2 + (1-\ga)^2 = 1 -2\ga +2\ga^2 = \frac12 +\frac12(1-2\ga)^2
\]
If $\ga<1/2$, then the diagonal entries of $M(t)$ are bounded away from zero.
Graphs with the above parameters exist if $n$ is a prime power.

\section{Questions?}

We do not know whether, for conference graphs, the diagonal entries of $M(t)$ can be 
bounded away from zero.

It might also be interesting to extend our results to walk-regular graphs (which can be 
characterized by the condition that the diagonals of the spectral idempotents are constant).

\bibliographystyle{plain}

\end{document}